\documentclass[12pt]{article}
\usepackage{a4,amsmath,amsthm,amssymb,amsfonts}

\theoremstyle{plain}

\newtheorem{thm}{Theorem}      

\newtheorem{lem}{Lemma}  
\newtheorem{prop}{Proposition} 

\theoremstyle{definition}
\newtheorem{prob}{Problem}

\newcommand{\C}{{\mathbb C}}
\newcommand{\abs}[1]{{\left| {#1} \right|}}
\newcommand{\p}[1]{{\left( {#1} \right)}}

\author{Johan Andersson\thanks{Department of Mathematics, Stockholm University, SE-10691, Sweden. {\it johana@math.su.se}}}

\begin{document}

\title{Explicit solutions to certain inf max problems from Tur{\'a}n power sum theory}


\maketitle

\begin{abstract}
 In a previous paper \cite{Andersson} we proved that
\begin{gather*}
 \sqrt n \leq \inf_{\abs{z_k} \geq 1}
 \max_{\nu=1,\ldots,n^2} \abs{\sum_{k=1}^n z_k^\nu} \leq \sqrt{n+1}
\\ \intertext{when $n+1$ is prime. In this paper we prove that}
 \inf_{ |z_k| = 1}
     \max_{\nu=1,\dots,n^2-n}   \left| \sum_{k=1}^n z_k^\nu \right| = \sqrt{n-1} \\ \intertext{when $n-1$ is a prime power, and}
 \inf_{\abs{z_k} \geq 1} \max_{\nu=1,\ldots,n^2-i} \abs{\sum_{k=1}^n z_k^\nu} =
  \sqrt n \qquad (i=2,\ldots,n-1) 
\end{gather*}
  when $n\geq 3$ is  a prime power.
 We give  explicit constructions
 of  $n$-tuples $(z_1,\ldots,z_n)$ which we prove are global minima for
  these problems.
  These are two of the few times in Tur\'an power sum theory where  solutions  in the $\inf \max$ problem can be explicitly
  calculated.
\end{abstract}

\section{Introduction}
In his paper  \cite{Turan2}  Tur\'an  shows that  
\begin{gather*}
 \inf_{\abs{z_k} \geq 1} \max_{\nu=1,\ldots,n} \abs{\sum_{k=1}^n z_k^\nu} =1.
\end{gather*}
Furthermore he gives an explicit construction
\begin{gather*}
  z_k=e \p{\frac {k} {n+1}} \qquad \qquad \p{e(x)=e^{2 \pi i x}}
\end{gather*}
which yields  a global minimum for the problem. He also showed that this is essentially (up to rearrangements of the $z_k$ and multiplication with a fixed unimodular constant) the only global minimum. Another result given by Cassels \cite {Cassels} is
\begin{gather}  \label{yt} 
  \inf_{z_1=1, \, z_k \in \C} \, \max_{\nu=1,\ldots,2n-1} \abs{\sum_{k=1}^n z_k^\nu} = 1.
\end{gather}
From this result it is clear that a global minimum can be given by $z_k=0$, for $k=2,\ldots,n$. In this case there are a lot of different solutions though. It has been shown by Dancs \cite{Dancs} that if we replace $2n-1$ by $3n-3$ there only exist the trivial solution. If we instead replace $2n-1$ by $3n-4$ there are already infinitely many different solutions to this problem.

Another problem studied in Tur\'an's book \cite{Turan} is 
\begin{gather*} 
   \inf_{z_1=1,\, z_k \in \C} \, \max_{\nu=1,\ldots,n} \abs{\sum_{k=1}^n z_k^\nu}.
\end{gather*} 
Atkinson \cite{Atkinson} showed that this quantity lies  in the interval $[1/6,1]$, Bir{\'o} (\cite{Biro1} and \cite{Biro2}) showed that it lies in the interval $[1/2+q,5/6]$ for some $q>0$ and  a sufficiently large $n$. In this case there exist no simple solution to the inf max problem, and in fact Cheer-Goldston \cite{CheerGoldston} used a computer to numerically obtain the minimal systems for small values of $n$. This is a typical situation in Tur\'an power sum theory. For most problems we have to be satisfied with inequalities and have little hope of obtaining an equality.

 In his book \cite{Turan} Tur\'an had a number of open problems. As problem 10 he proposed the study of the quantity
\begin{gather} \label{iii}
   \inf_{\abs{z_k} \geq 1}\max_{\nu=1,\ldots,n^2} \abs{\sum_{k=1}^n z_k^\nu}.
\end{gather}
In a previous paper \cite{Andersson} we proved the strong inequality
\begin{gather} \label{ii}
 \sqrt n \leq \inf_{\abs{z_k} \geq 1} \max_{\nu=1,\ldots,n^2} \abs{\sum_{k=1}^n z_k^\nu} \leq \sqrt{n+1}
\end{gather}
which is valid whenever $n+1$ is prime. The proof relied on an explicit construction due to Hugh Montgomery (see Tur\'an \cite{Turan} p. 83).  It would be interesting to find a true global minimum for the problem, and we may ask:
\begin{prob} \label{prob1}
    Does Montgomery's construction (for some sufficiently large $n$ such that $n+1$ is prime) give a true global minimum for the quantity \eqref{iii}, or equivalently can the upper inequality be replaced by an equality in eq. \eqref{ii}?
\end{prob}
We will not be able to answer this question in this paper. We will however be able to answer the corresponding question in some closely related problems. We will use explicit constructions to prove
  the following theorems:

\begin{thm} \label{ajj} Suppose that $n-1$ is a prime power. Then
 \begin{gather*}
 \inf_{ |z_k| = 1}
     \max_{\nu=1,\dots,n^2-n}   \left| \sum_{k=1}^n z_k^\nu \right| = \sqrt{n-1}.
\end{gather*}
Furthermore an explicit $n$-tuple $(z_1,\ldots,z_n)$ which gives a global minimum is given by Theorem \ref{Singer} and eq. \eqref{yyy3}.
\end{thm}

\begin{thm} \label{and} Let $n \geq 3$ be a prime power, and let $2 \leq i \leq n-1$. Then
\begin{gather*}
  \inf_{\abs{z_k} \geq 1} \max_{\nu=1,\ldots,n^2-i} \abs{\sum_{k=1}^n z_k^\nu} = \sqrt
  n.
\end{gather*}
Furthermore the explicit $n$-tuple $(z_1,\ldots,z_n)$ given by Theorem \ref{Bose}
 and eq. \eqref{yyy} provides a global minimum for the problem.
\end{thm}

\section{Fabykowski's construction}

In our paper \cite{Andersson} we also used a construction of Fabrykowski \cite{Fabrykowski} to prove some results closely related to eq. \eqref{ii}. The construction by Fabrykowski
depends of the existence of a perfect difference set
\begin{thm} {\rm (Singer)} \label{Singer} Let $n-1$ be a prime power. Then there exists integers $a_1,\ldots,a_{n}$
such that the integers $a_i-a_j$ for $i \neq j$ form all non zero
residues $\mod n^2-n+1$.
\end{thm}
 \noindent of Singer \cite{Singer}. 
\noindent By choosing \begin{gather} \label{yyy3} z_k=e \p{\frac {a_k}{n^2-n+1}}, \qquad \p{e(x)=e^{2 \pi i x}} \\
\intertext{we see that} \notag \begin{split} \abs{\sum_{k=1}^n z_k^\nu}^2 &=
n+\sum_{i \neq j} e \p{\frac {\nu(a_j-a_i)}{n^2-n+1}}  \\ &= n-1 +
\sum_{j=1}^{n^2-n+1} e \p{\frac {\nu j}{n^2-n+1}} \\ &= \begin{cases} n-1, &  n^2-n+1 \not |\nu,
\\ n^2, & n^2-n+1 |\nu. \end{cases}
\end{split}
\end{gather}
Hence we obtain (Andersson \cite{Andersson} Lemma 3) 
\begin{lem} \label{fab} Let $n$ be a prime power. There exists
an n-tuple of unimodular complex numbers such that 
\begin{gather*}
 \abs{\sum_{k=1}^n z_k^\nu} =\begin{cases} \sqrt {n-1}, & n^2-n+1 \not | \nu,
 \\ n, & n^2-n+1 | \nu. \end{cases}
\end{gather*}
\end{lem}
In Andersson \cite{Andersson} we used a lemma of Cassels to prove that (\cite{Andersson}, Corollary 1)
\begin{gather} \label{hhhu}
  \inf_{|z_k| \geq 1} \max_{\nu=1,\ldots,2nm-m(m+1)+1} \abs{\sum_{k=1}^n z_k^\nu} \geq \sqrt{m}. \qquad (1 \leq m  \leq n)
\end{gather}
Together with Lemma \ref{fab}  this implies (\cite{Andersson}, Proposition 1 $(ii)$) 
\begin{gather}\label{hhh}
   \sqrt{n-2} \leq \inf_{ |z_k| \geq 1}
     \max_{\nu=1,\dots,n^2-n} 
     \left| \sum_{k=1}^n z_k^\nu \right| \leq \sqrt{n-1}.
 \end{gather}
There exist another way of getting lower bounds for the inf max problem, which is a theorem independently proved by Newman, Cassels and Szalay \cite{Szalay}
 (see Theorem 7.3 in Tur\'an \cite{Turan}). The result is more general than eq. \eqref{hhhu} since it does not assume that we consider the  pure power sum problem (We can also have coefficients $b_j>0$). However it is also less general since it assumes that $|z_k|=1$. We will state the result for the pure power sum case:  
\begin{lem} Suppose that $z_k$ are unimodular complex numbers, and $c\geq 1$ is an integer. Then
\begin{gather*}
 \max_{1 \leq \nu \leq c n} \left| \sum_{k=1}^n z_k^\nu \right| \geq \sqrt{\frac {cn-n+ 1} {c}}.
\end{gather*}
\end{lem}
\noindent In the special case $c=n-1$ we obtain the lower bound
\begin{gather*}
 \max_{1 \leq \nu \leq n^2-n} \abs{\sum_{k=1}^n z_k^\nu} \geq  \sqrt {n-1}
\end{gather*}
\noindent and by combining this with Lemma \ref{fab} we obtain a proof of Theorem \ref{ajj}. 
Natural problems to ask  are:
\begin{prob} \label{prob2}
 Does all global minima of the min max problem in Theorem \ref{ajj} arise
 from Fabrykowski's  construction (and multiplication of a fixed unimodular constant)?
 \end{prob}
\begin{prob} \label{prob7}
 Can the condition $|z_k|=1$ in Theorem \ref{ajj} be replaced by $|z_k| \geq 1$ ? 
\end{prob}
By the result of  Blanksby \cite{Blanksby} (or it can be proved immediately from Tur\'an's second main theorem \cite{Turan} Chapter 8) we can give a partial answer to Problem \ref{prob7}:
\begin{prop}
 There exist a constant $C>0$ such that for every $n \geq 1$ we have that if $(z_1,\ldots,z_n)$ is a global minimum for the problem in eq. \eqref{hhh}, then $|z_k| \leq (1+C/n)$.
\end{prop}

\section{A new construction and a theorem of Bose}

This  theorem of Singer is also used to
construct  Golomb rulers and Sidon sets (For a survey, see Dimitromanolakis \cite{dimitromanolakis}), and in fact
Montgomery's construction  is in this setting equivalent
to Ruzsa's construction (see \cite{Ruzsa} Theorem 2) which is also
used to construct Golomb rulers and Sidon sets. There also exists
 a third construction of Bose \cite{Bose} which is used to construct
Golomb rulers that has hitherto not been used in the corresponding
power sum problem. We will now  see what this
construction will yield when applied to the power sum problem.

We first state a result taken from Bose \cite{Bose}:
\begin{thm}{\rm (Bose)} \label{Bose} Let $n$ be a prime power. There exists integers $b_1,\ldots,b_n$
such that the residues $b_i-b_j$ for $i \neq j$ form all  residues
 $\mod n^2-1$ which are not divisible by $n+1$.
\end{thm}

\noindent By choosing \begin{gather} \label{yyy} z_k=e \p{\frac {b_k}{n^2-1}}, \qquad \p{e(x)=e^{2 \pi i x}} \\
\intertext{we see that} \notag \begin{split} \abs{\sum_{k=1}^n z_k^\nu}^2 &=
n+\sum_{i \neq j} e \p{\frac {\nu(b_j-b_i)}{n^2-1}}  \\ &= n+
\sum_{j=1}^{n^2-1} e \p{\frac {\nu j}{n^2-1}}- \sum_{j=1}^{n-1} e
\p{\frac {\nu j}{n-1}} \\ &= \begin{cases} n, &  (n-1) \not |\nu,
\\ 1, & (n-1) |\nu. \end{cases} \qquad (\nu=1, \dots,n^2-2)
\end{split}
\end{gather}
Hence we obtain \begin{lem} Let $n$ be a prime power. There exists
an n-tuple of unimodular complex numbers such that for
$\nu=1,\dots,n^2-2$  one has that
\begin{gather}
 \abs{\sum_{k=1}^n z_k^\nu} =\begin{cases} \sqrt n, & (n-1) \not | \nu,
 \\ 1, & (n-1)| \nu. \end{cases}
\end{gather}
\end{lem}
\noindent By combining  this with the choice $m=n$ in eq. \eqref{hhhu} (see also Andersson \cite{Andersson}, Corollary 3)  we obtain a proof of Theorem \ref{and}.
As in Problem \ref{prob2} we may ask
\begin{prob} 
 Does all global minima of the min max problem in Theorem \ref{and} arise
 from Bose's  construction (and multiplication of a fixed unimodular constant)?
 \end{prob}
 All problems 1-4 would be interesting to investigate numerically as in
 Cheer-Goldston \cite{CheerGoldston} for small values of $n$.

\section{Acknowledgments}

I would like to thank  Alexei Venkov for inviting me to Aarhus
and asking me some questions which got me thinking  about this problem again. I would also like to thank the referee for some suggestions on how to improve the exposition.

\bibliographystyle{plain}

\begin{thebibliography}{10}

\providecommand{\bysame}{\leavevmode\hbox to3em{\hrulefill}\thinspace}

\bibitem{Andersson}
J.~Andersson.
\newblock On some power sum problems of {T}ur\'an and {E}rd{\H o}s.
\newblock {\em Acta Math. Hungar.}, 70(4):305--316, 1996.

\bibitem{Atkinson}
F.~V. Atkinson.
\newblock On sums of powers of complex numbers.
\newblock {\em Acta Math. Acad. Sci. Hungar.}, 12:185--188, 1961.

\bibitem{Biro2}
A.~Bir{\'o}.
\newblock An upper estimate in {T}ur\'an's pure power sum problem.
\newblock {\em Indag. Math. (N.S.)}, 11(4):499--508, 2000.

\bibitem{Biro1}
\bysame
\newblock An improved estimate in a power sum problem of {T}ur\'an.
\newblock {\em Indag. Math. (N.S.)}, 11(3):343--358, 2000.

\bibitem{Blanksby}
P.~E. Blanksby.
\newblock Sums of powers of conjugates of algebraic numbers.
\newblock {\em Proc. Amer. Math. Soc.}, 49:28--32, 1975.

\bibitem{Bose}
R.~C. Bose.
\newblock An affine analogue of {S}inger's theorem.
\newblock {\em J. Indian Math. Soc. (N.S.)}, 6:1--15, 1942.

\bibitem{Cassels}
J.~W.~S. Cassels.
\newblock On the sums of powers of complex numbers.
\newblock {\em Acta Math. Acad. Sci. Hungar.}, 7:283--289, 1956.

\bibitem{CheerGoldston}
A.~Y. Cheer and D.~A. Goldston.
\newblock Tur\'an's pure power sum problem.
\newblock {\em Math. Comp.}, 65(215):1349--1358, 1996.

\bibitem{Dancs}
I.~Dancs.
\newblock Power sums of complex numbers.
\newblock {\em Mat. Lapok}, 13:108--114, 1962.

\bibitem{dimitromanolakis}
A.~Dimitromanolakis.
\newblock Analysis of the {G}olomb ruler and the {S}idon set problems and
  determination of large near-optimal golomb rulers.
\newblock http://citeseer.ist.psu.edu/dimitromanolakis02analysis.html

\bibitem{Fabrykowski}
J.~Fabrykowski.
\newblock A note on sums of powers of complex numbers.
\newblock {\em Acta Math. Hungar.}, 62(3-4):209--210, 1993.

\bibitem{Ruzsa}
I~Z. Ruzsa.
\newblock Solving a linear equation in a set of integers. {I}.
\newblock {\em Acta Arith.}, 65(3):259--282, 1993.

\bibitem{Singer}
J.~Singer.
\newblock A theorem in finite projective geometry and some applications to
  number theory.
\newblock {\em Trans. Amer. Math. Soc.}, 43(3):377--385, 1938.

\bibitem{Szalay}
M.~Szalay.
\newblock {\em On number theoretical extremal problems}.
\newblock Thesis, 1974.

\bibitem{Turan2}
P.~Tur{\'a}n.
\newblock On a certain limitation of eigenvalues of matrices.
\newblock {\em Aequationes Math.}, 2:184--189, 1969.

\bibitem{Turan}
\bysame
\newblock {\em On a new method of analysis and its applications}.
\newblock Pure and Applied Mathematics (New York). John Wiley \& Sons Inc., New
  York, 1984.
\newblock With the assistance of G. Hal\'asz and J. Pintz, With a foreword by
  Vera T. S\'os, A Wiley-Interscience Publication.

\end{thebibliography}

\end{document}